\documentclass[11pt]{article}
\usepackage[english]{babel}    
\usepackage[utf8]{inputenc}    
\usepackage{blindtext}         
\usepackage[T1]{fontenc}       

\usepackage[sumlimits,%
                intlimits,%
                namelimits%
                ]{amsmath}         
\usepackage{amsthm}            
\usepackage{amscd}             

\usepackage{fouridx}
\usepackage{mathtools}

\usepackage{amssymb}

\usepackage[hyperfootnotes=false]{hyperref}
\usepackage{cleveref}
\usepackage{autonum}
\usepackage{enumerate}

\newtheorem{theorem}{Theorem}

\def\RR{{\mathbb{R}}}

\def\phi{\varphi}
  
\begin{document}
\begin{center}
\Large
 \bf High Accuracy Quasi-Interpolation using a new class of generalized Multiquadrics 
\end{center}
\bigskip

\begin{center}
 Dedicated to Dany Leviatan on the occasion of his 80th birthday\\ \bigskip \bigskip
 \it Mathis Ortmann, Justus-Liebig University,
  Mathematics, 35392 Giessen, and \\ \smallskip
  Martin Buhmann, Justus-Liebig University, Mathematics, 35392 Giessen

\end{center}
\bigskip\bigskip
{\small {\it Abstract: \/} A new generalization of multiquadric functions $\phi(x)=\sqrt{c^{2d}+||x||^{2d}}$, where $x\in\mathbb{R}^n$, $c\in \mathbb{R}$, $d\in \mathbb{N}$, is presented to increase the accuracy of quasi-interpolation further. With the restriction to Euclidean spaces of odd dimensionality, the generalization can be used to generate a quasi-Lagrange operator that reproduces all polynomials of degree $2d-1$. 
In contrast to the classical multiquadric, the convergence rate of the quasi-interpolation operator can be significantly improved by a factor $h^{2d-n-1}$, where $h>0$ represents the grid spacing. Among other things, we compute the generalized Fourier transform of this new multiquadric function. Finally, an infinite regular grid is employed to analyse the properties of the aforementioned generalization in detail.}

\date{\today}  
 \bigskip\bigskip\bigskip
 \title{Introduction}
\section{Introduction}

Quasi-interpolation with radial basis functions (RBFs) is a powerful technique widely used in numerical analysis and approximation theory. It provides a flexible and efficient approach for constructing approximate functions based on given data points without explicitly performing interpolation. Applications of quasi-interpolation with RBFs include surface reconstruction, data fitting, image processing, and solving partial differential equations \cite{WOS:000305709900003},\cite{WOS:000513922200014},\cite{WOS:000442220800001}. It finds utility in various fields, such as computer graphics, computational physics, geostatistics and machine learning, \cite{WOS:000319478400017}, \cite{WOS:000479317200007}.
Unlike traditional interpolation methods that aim to find a global polynomial, spline or RBF-interpolant $s(x)=\sum_{\xi \in \Xi} a_\xi \phi(||x-\xi||)$ that passes through all data points $\Xi$, quasi-interpolation with RBFs focus on local Taylor approximations in order to approximate differential functions \cite{qi}. The key idea is to express the approximate function as a linear combination of a \textit{quasi-Lagrange-function} at specific points. These Lagrange functions possess desirable properties, giving the approximants polynomial reproduction, fast decay and smoothness. 
One way to construct such quasi-Lagrange-functions is the RBF-approach, which is the focus of this paper. We will construct quasi-Lagrange functions of the form
\begin{align}
\Psi_\alpha(x)=\sum_{i\in\mathbb{Z}^n}\mu_{\alpha,i}\phi(||x-i||),
\end{align}
where $x\in\mathbb{R}^n,\alpha\in \mathbb{A}$, where $\mathbb{A}$ is some countable subset of $\mathbb{R}^n$, $\phi\in \mathrm{C}^\infty(\mathbb{R})$ is a RBF and $||\cdot||$ is the Euclidean norm. The sum over $i$ can be finite or infinite, depending on the dimension of $x$ and the choice of RBFs \cite{qi}. The approximation to given data points $\mathbb{A} \subset \mathbb{R}^n$ of a function $f$ is then given as $Qf(x)=\sum_{\alpha \in\mathbb{A}}f(\alpha)\Psi_{\alpha}(x)$.
In order to guarantee its applicability at a minimum we require that $\Psi_\alpha (x)$ form a partition of unity. The choice of RBFs, such as the Gaussian function or the multiquadric function, and their associated weights $\mu_{\alpha,i}$ are crucial in quasi-interpolation. The most often used RBF is the multiquadric but it is only capable of reproducing polynomials of degree $n$ in $n$ dimensions \cite{BUHMANN1988}.
Therefore we introduce a new generalization of multiquadrics, namely
\begin{align}
\phi(x)=\sqrt{c^{2d}+||x||^{2d}},\label{eq: genRBF}
\end{align} 
where $c\in\mathbb{R}$ and $d \in \mathbb{N}$. We believe that this RBF has not been considered before. Given this new RBF, we will construct a quasi-Lagrange function with sufficient decay rate such that the quasi-interpolant reproduces polynomials of degree $2d-1$. Given this information, the famous \textit{ Strang and Fix Conditions} given below, will lead to the approximation order. To keep the results more comprehensible and manageable we will perform all calculations on a regularly spaced grid in $\mathbb{R}^n$ with spacing $h$. For a multi-index $\alpha$, $D^\alpha=\frac{\partial^{\alpha_1}}{\partial x_1^{\alpha_1}} \cdots \frac{\partial^{\alpha_m}}{\partial x_m^{\alpha_m}}$ will denote the higher order partial derivatives in the following.

\begin{theorem}\label{SF}[Strang and Fix conditions]

Let $m$ be a positive integer and $\Psi:\mathbb{R}^n \rightarrow \mathbb{R}$ be a function such that
\begin{enumerate}
\item there exists a non-negative real valued $\ell$ such that, when $\|x\|\rightarrow \infty$, $|\Psi(x)|=\mathcal{O}(\Vert x\Vert^{-n-m-\ell})$ and this implies $\hat\Psi \in C^{m+\ell-1}(\mathbb{R}^n),$
\item $D^{\alpha} \hat{\Psi}(0)=0$, $\forall \alpha \in \mathbb{Z}_+^n$, $1\leq \vert \alpha\vert\leq m$, and $\hat{\Psi}(0)=1$,
\item $ D^{\alpha}\hat{\Psi}(2\pi j)=0,\ \forall j \in \mathbb{Z}^n\setminus\{0\}$ and  $\forall \alpha \in \mathbb{Z}_+^n$ with $ \vert\alpha\vert \leq m$. 
\end{enumerate}

Then the quasi-interpolant
\begin{equation}
Q_hf(x)=\sum_{j\in\mathbb{Z}^n}f(jh)\Psi(x/h-j),\qquad x\in\RR^n,
\end{equation}
is well-defined and exact on the space of polynomials of degree $m$ and the uniform approximation error can be estimated by
$$\Vert Q_hf-f\Vert_{\infty}=\begin{cases}\mathcal{O}(h^{m+\ell}), & \text{when } 0<\ell<1,\\
\mathcal{O}(h^{m+1}\log (1/h)), & \text{when } \ell=1,\\
\mathcal{O}(h^{m+1}), & \text{when } \ell>1,\\
\end{cases}$$
for $h\rightarrow 0$ and a bounded function $f\in C^{m+1}(\RR^n)$ with bounded derivatives \cite{BUHMANN2015156}.
\end{theorem}

\section{Fourier transform}
To employ the \textit{Strang and Fix conditions}, one needs the Fourier transform of the RBF. The $n$-dimensional Fourier transform is defined as
\begin{align}
\hat f(x):=\int_{\mathbb{R}^n}f(y)\textnormal{e}^{-\mathrm{i}x\cdot y} \mathrm{d}y ,\qquad x\in\mathbb{R}^n \, ,
\end{align}
where $f$ is an integrable function. Since $\phi$ in equation~(\ref{eq: genRBF}) is not integrable, we need to employ the theory of generalized Fourier transforms \cite{jones}.

\begin{theorem}\label{FT}
The generalized $n$-dimensional Fourier transform of $$\phi(x)=\sqrt{c^{2d}+||x||^{2d}}\qquad  c\in\mathbb{R},\, d\in \mathbb{N}, \,x \in \mathbb{R}^n$$ is given by
\begin{align}
\hat\phi(s)=- 2^{n-1} \pi ^{\frac{n-1}{2}} d^{n/2} c^d s^{-n}G_{0,2d}^{\,d+1,0}\!\left(\left.{\begin{matrix}-\\ -\frac{1}{2},\frac{n}{2d},\frac{n+2}{2d},\dots \frac{2d-2+n}{2d},\frac{1}{d},\frac{2}{d}\dots \frac{d-1}{d}\end{matrix}}\;\right|\,\left(\frac{cs}{2d}\right )^{2d}\right),
\end{align}
where $s\in\mathbb{R}_+=\{x\in\mathbb{R}|\,x>0\}$ is the radial part of its argument and $G$ is the Meijer G-function. Alternatively, the Meijer G-function can be generalized to a Fox H-function. Then the Fourier transform is given by
\begin{align}
\hat\phi(s)=- 2^{n-1} \pi ^{\frac{n-1}{2}} d^{n/2} c^d s^{-n}H_{1,0}^{\,2,1}\!\left(\left.{\begin{matrix}\left(1,1\right)\\ \left(-\frac{1}{2},1\right),\left(-\frac{n}{2},d\right)\left(1,d\right)\end{matrix}}\;\right|\,\left(\frac{cs}{2}\right )^{2d}\right).
\end{align}
The definition of the Meijer G- and Fox H-function are given in the proof.
\end{theorem}
\subsection*{\textit{Note}}
In the special case $d=1$ the Meijer G-function reduces to the modified Bessel function of the second kind, and one obtains the well known generalized Fourier transform of the multiquadric. Also the case $c=0$ is allowed but provides the classical results only. Also the Fox H-function representation holds true for $d \in \mathbb{R}_+$, while die Meijer G-function representation is only valid for $d \in \mathbb{N}$. Still, these functions are merely names for the particular Fourier transform and the proof will rely on its integral representation.

\subsection*{\textit{Proof}}
The Fourier transform of a radial symmetric function is also a radial symmetric function. Thus we can write

\begin{align}
\hat\phi (s)= \lim_{\varepsilon \to 0_+}\frac{(2\pi)^{\frac{n}{2}}} {s^{\frac{n}{2}-1}}\int_0^\infty r^\frac{n}{2}\phi(r) J_{\frac{n}{2}-1}(sr) \mathrm{e}^{-\varepsilon r^2} \, \mathrm{d}r, \qquad s\in \mathbb{R}_+ ,
\end{align}
where we introduced the Gauss function as a convergence generating factor to keep the integral finite. 
Normally this refers to the Hankel transform of order $\frac{n}{2}-1$, but here it is more appropriate to think of the integral as a Mellin transform.
Rewriting the RBF and the Bessel function in terms of Meijer G-functions yields

\begin{align}
\phi(r)=\sqrt{c^{2d}+r^{2d}}&=-\frac{c^d}{2\sqrt{\pi}}G_{1,1}^{\,1,1}\!\left(\left.{\begin{matrix}\frac{3}{2}\\ 0\end{matrix}}\;\right|\,\left(\frac{r}{c}\right )^{2d}\right),\\
J_{\frac{n}{2}-1}(sr)&=G_{0,2}^{\,1,0}\!\left(\left.{\begin{matrix}-\\ \frac{n-2}{4},-\frac{n-2}{4}\end{matrix}}\;\right|\,\frac{(sr)^2}{4}\right), \, r>0.
\end{align}
The Meijer G-function is defined as
\begin{align}
G_{p,q}^{\,m,n}\!\left(\left.{\begin{matrix}a_{1},\dots ,a_{p}\\b_{1},\dots ,b_{q}\end{matrix}}\;\right|\,z\right)={\frac {1}{2\pi \mathrm{i}}}\int _{L}{\frac {\prod _{j=1}^{m}\Gamma (b_{j}-s)\prod _{j=1}^{n}\Gamma (1-a_{j}+s)}{\prod _{j=m+1}^{q}\Gamma (1-b_{j}+s)\prod _{j=n+1}^{p}\Gamma (a_{j}-s)}}\,z^{s}\,\mathrm{d}s,
\end{align}
where the path of integration can have three different shapes \cite{handbook}. In our case, $L$ is a loop that starts at infinity on a line parallel to the positive real axis, encircles the poles of the $\Gamma (b_\ell-s)$ once in the negative sense and returns to infinity on another line parallel to the positive real axis \cite{olv}.
The Meijer G-function is a Mellin-Barnes type integral and can be viewed as an inverse Mellin transform \cite{Carlos}, because
\begin{align}
\int _{0}^{\infty }z^{s-1}\;G_{p,q}^{\,m,n}\!\left(\left.{\begin{matrix}a_{1},\dots ,a_{p}\\b_{1},\dots ,b_{q} \end{matrix}}\;\right|\,\eta z\right)\mathrm{d}z={\frac {\eta ^{-s}\prod _{j=1}^{m}\Gamma (b_{j}+s)\prod _{j=1}^{n}\Gamma (1-a_{j}-s)}{\prod _{j=m+1}^{q}\Gamma (1-b_{j}-s)\prod _{j=n+1}^{p}\Gamma (a_{j}+s)}},
\end{align}
where $s,\eta \in\mathbb{R}$.
Substituting $\tilde r =r^2$ and using the identity
\begin{align}
z^{\rho }\;G_{p,q}^{\,m,n}\!\left(\left.{\begin{matrix}\mathbf {a_{p}} \\\mathbf {b_{q}} \end{matrix}}\;\right|\,z\right)=G_{p,q}^{\,m,n}\!\left(\left.{\begin{matrix}\mathbf {a_{p}} +\rho \\\mathbf {b_{q}} +\rho \end{matrix}}\;\right|\,z\right)
\end{align}
on the Bessel functions yields

\begin{align}
\hat\phi (s)=-  \frac{2^{n-1}\pi^{\frac{n-1}{2}}c^d} {4s^{n-2}}\lim_{\varepsilon \to 0_+}\int_0^\infty  \mathrm{e}^{-\varepsilon \tilde r} G_{1,1}^{\,1,1}\!\left(\left.{\begin{matrix}\frac{3}{2}\\ 0\end{matrix}}\;\right|\,\frac{\tilde r^d}{c^{2d}}\right) G_{0,2}^{\,1,0}\!\left(\left.{\begin{matrix}-\\ \frac{n}{2}-1,0\end{matrix}}\;\right|\,\frac{s^2\tilde r}{4}\right) \mathrm{d}\tilde r.
\end{align}
Writing the first Meijer G-function as a Mellin-Barnes integral, changing the order of integration results in
\begin{align}
\hat\phi (s)=-  \frac{2^{n-1}\pi^{\frac{n-1}{2}}c^d} {4s^{n-2}}\lim_{\varepsilon \to 0_+} \frac{1}{2\pi \mathrm{i}} \int_L\Gamma \left(-t\right) \Gamma \left(-\frac{1}{2}+t\right)  
\times{}\\{}\times \int_0^\infty\left(\frac{\tilde r}{c^2}\right)^{dt}  \mathrm{e}^{-\varepsilon \tilde r}
G_{0,2}^{\,1,0}\!\left(\left.{\begin{matrix}-\\ \frac{n}{2}-1,0\end{matrix}}\;\right|\,\frac{s^2\tilde r}{4}\right) \mathrm{d}\tilde r \, \mathrm{d}t.
\end{align}
The change of integration is valid through Fubini's theorem.  Additionally, the inner integral can be found in \cite{Prudnikov}, providing
\begin{align}
\hat\phi (s)=-  \frac{2^{n-1}\pi^{\frac{n-1}{2}}c^d} {4s^{n-2}}\lim_{\varepsilon \to 0_+} \frac{1}{2\pi \mathrm{i}}\int_L\frac{\Gamma \left(-t\right) \Gamma \left(-\frac{1}{2}+t\right)}{c^{2dt}} \varepsilon^{-dt-1} G_{2,1}^{\,1,1}\!\left(\left.{\begin{matrix}-dt\\ \frac{n}{2}-1,0\end{matrix}}\;\right|\,\frac{s^2}{4\varepsilon}\right)  \mathrm{d}t.
\end{align}
Using the Beppo Levi monotone convergence theorem, the integral and the limit can be exchanged. To evaluate the limit, we  write the Meijer G-function as a Mellin-Barnes integral and expand it into a series using the residue theorem. It is easy to see that in the limit $\varepsilon \to 0_+$ the only nonvanishing term is induced by the first pole at $dt+1$. Ultimately, we obtain the integral representation.

\begin{align}
\hat\phi (s)=-  \frac{2^{n-1}\pi^{\frac{n-1}{2}}c^d} {s^n}\frac{1}{2\pi \mathrm{i}}\int_L\frac{ \Gamma \left(-t\right ) \Gamma \left(-\frac{1}{2}+t\right )  \Gamma \left(\frac{n}{2}+dt\right )}{ \Gamma \left(-dt\right )} \left(\frac{2}{sc}\right)^{2dt} \mathrm{d}t. \label{eq: FT}
\end{align}
Comparing with the definition of the Fox H-function
\begin{align}
 &H_{p,q}^{\,m,n}\!\left(\left.{\begin{matrix}\left(a_1,A_1\right),\dots,\left(a_p,A_p\right)\\ \left(b_1,B_1\right),\dots, \left(b_q,B_q\right)\end{matrix}}\;\right|\,z\right)\\
&={\frac {1}{2\pi \mathrm{i}}}\int _{L}{\frac { \prod _{j=1}^{m}\Gamma (b_{j}+B_{j}s)\,\prod _{j=1}^{n}\Gamma (1-a_{j}-A_{j}s)}{ \prod _{j=m+1}^{q}\Gamma (1-b_{j}-B_{j}s)\,\prod _{j=n+1}^{p}\Gamma (a_{j}+A_{j}s)}}z^{-s}\,ds,
\end{align}
where $A_1\dots A_p$ and $B_1\dots B_q$ are positive numbers, the desired representation of the Fox H-function is obtained.
To expand the integral representation into a Meijer G-function, we utilize the Gauss multiplication formula \cite{Abramowitz}
\begin{align}
\frac{(2\pi)^{(d-1)/2}}{d^{x-1/2}}\cdot \Gamma\left ( x\right)=\Gamma \left(\frac{x}{d}\right ) \cdot \Gamma \left(\frac{x+1}{d}\right ) \cdots \Gamma \left(\frac{x+d-1}{d}\right )
\end{align}
and apply
\begin{align}
G_{p,q}^{\,m,n}\!\left(\left.{\begin{matrix}\mathbf {a_{p}} \\\mathbf {b_{q}} \end{matrix}}\;\right|\,z\right)=G_{q,p}^{\,n,m}\!\left(\left.{\begin{matrix}1-\mathbf {b_{q}} \\1-\mathbf {a_{p}} \end{matrix}}\;\right|\,z^{-1}\right).
\end{align}
The Gauss multiplication formula is only applicable when $d$ is an integer. It is solely used to convert the integral representation from equation (\ref{eq: FT}) to a Meijer G-function. However, the integral- and the Fox H-function representation remains valid for $d \in \mathbb{R}_+$. For additional information on identities associated with Meijer G-functions, please refer to \cite{luke1975mathematical}.

\begin{theorem}
The function $\phi(x)=\sqrt{c^{2d}+||x||^{2d}}$ , $c\in\mathbb{R}$, $ d\in \mathbb{N}$, $x \in \mathbb{R}^n$, is not positive definite.
\end{theorem}

\subsection*{\textit{Proof}}
Let $x_1$ and $x_2$ be two distinct points in $\mathbb{R}^n$. Then the interpolation matrix is given by
$$A=
\begin{bmatrix}
c^d & \sqrt{c^{2d}+||x_1-x_2||^{2d}} \\
\sqrt{c^{2d}+||x_1-x_2||^{2d}} & c^d
\end{bmatrix}.$$
Furthermore the eigenvalues of the matrix $A$ are given by 
\begin{align}
\lambda_1 &=\sqrt{c^{2 d}}- \sqrt{||x_1-x_2||^{2d}+c^{2 d}} <0\\
\lambda_2 &=\sqrt{c^{2 d}}+\sqrt{||x_1-x_2||^{2d}+c^{2 d}} >0.
\end{align}
From this it follows that $A$ is not positive definite and the same applies to $\phi(x)$.


\section{Asymptotic behaviour of $\hat\phi(s)$}
Next we analyse the asymptotic behaviour of $\hat\phi(s)$ as $s \to 0_+$. 

\begin{theorem}\label{Satz: asymptoticbehaviour}
Let $\hat\phi(s)$ be the generalized Fourier transform of $\phi(x)=\sqrt{c^{2d}+||x||^{2d}}$, where $n$ is the dimension and $d \in \mathbb{R}_+$ is the generalization parameter. Then

\begin{align}
\hat\phi(s)= 
\begin{cases}
\mathcal{O}\left(1\right)\qquad &\textnormal{as } s\to 0 \textnormal{ for } d \textnormal{ even}  \textnormal{ and }n< d\\
\mathcal{O}\left(s^{-n+d}\right) \qquad &\textnormal{as } s\to 0 \textnormal{ for } d \textnormal{ even}  \textnormal{ and }n >d\\
\mathcal{O}\left(-\log (s)\right)\qquad &\textnormal{as } s\to 0 \textnormal{ for } d \textnormal{ even}  \textnormal{ and }n=d\\
\mathcal{O}\left(s^{-n-d}\right)\qquad &\textnormal{as } s\to 0 \textnormal{ otherwise}.
\end{cases}
\end{align}
The exact asymptotic behaviour for $s \to 0_+$ can be found in equations ~(\ref{eq: d_not_even}),(\ref{eq: d_even_nkd}),(\ref{eq: d_even_ngd}) and (\ref{eq: d_even_ngleichd}).
\end{theorem}

\subsection*{\textit{Note}}
The Meijer G representation of the Fourier transform is only valid if $d$ is an integer. The results given are valid for $d \in \mathbb{R}_+$.

\subsection*{\textit{Proof}}
The integral representation of the Fourier transform is given in equation~(\ref{eq: FT}). To analyse the asymptotic behaviour, we consider only the integral part. Therefore we define
\begin{align}
G:=\frac{1}{2\pi \mathrm{i}}\int_L\frac{\Gamma \left ( - \frac{1}{2}-t \right)\Gamma \left (t \right)\Gamma \left (  \frac{n}{2}-dt \right)}{\Gamma \left ( dt \right)}\left( \frac{cs}{2}\right)^{2dt} \mathrm{d}t\, , \label{eq: G}
\end{align}
where $L$ is the path coming from $+\infty-i\varepsilon$ to $-\frac{1}{2}-i\varepsilon$, circles the pole at $-\frac{1}{2}$ and going back to $+\infty +i\varepsilon$. 
This is the most important formula for determining the asymptotic behaviour of our Fourier transform.
The asymptotic behaviour of $G$ for $s \to 0$ is given by the residue of the integral at the pole $t=-1/2$. Keep in mind that the integration path circles the poles in a negative sense, such that a minus sign appears when the residue theorem is used.
Simple poles can easily evaluated using the identity
\begin{align}
\textnormal{Res}(\Gamma ,-z)=\frac{(-1)^z}{z!}\, ,
\end{align}
where $z \in \mathbb{N}$. We follow the asymptotic notation used in \cite{Bender.1999} and define the asymptotic symbol "$\sim$" as
\begin{align}
f(x)\sim g(x)\textnormal{ as } x\to x_0 \Leftrightarrow \lim_{x\to x_0} \frac{f(x)}{g(x)} =1.
\end{align}

\subsubsection*{Case I: $d$  not an even integer}
Let $d$ be not an even integer. Then the asymptotic behaviour of $G$ and $\hat \varphi(s)$ is given by
\begin{align}
G\sim \frac{\sqrt{\pi } 2^{d+1}   \Gamma \left(\frac{d+n}{2}\right)}{\Gamma \left(-\frac{d}{2}\right)} (c s)^{-d}\qquad \textnormal{as } s\to 0_+
\end{align}
and
\begin{align}
\hat\phi (s)\sim -2^{n+d} \pi ^{\frac{n}{2}} \frac{\Gamma \left(\frac{d+n}{2}\right)}{\Gamma \left(-\frac{d}{2}\right)} s^{-n-d} =\mathcal{O}(s^{-n-d})\qquad \textnormal{as } s\to 0_+\, . \label{eq: d_not_even}
\end{align}
 If $d$ is even the denominator becomes singular and the term vanishes.

\subsubsection*{Case II: $d$ is an even integer}
 For even $d$ the first nonvanishing pole is $ t=\min\left\{\frac{1}{2},\frac{n}{2d}\right\}$. One has to be careful to see whether this pole is simple or double. 
If $ n<d$ then the pole is simple and we obtain the asymptotic behaviour of $G$ and $\hat\varphi(s)$ by

\begin{align}
G\sim\frac{\Gamma \left(-\frac{1}{2}-\frac{n}{2d}\right)\Gamma \left(\frac{n}{2d}\right)}{d \Gamma \left(\frac{n}{2}\right)}\left( \frac{cs}{2}\right)^{n}\qquad \textnormal{as } s\to 0_+ 
\end{align}
and 
\begin{align}
\hat\phi (s)\sim-\frac{ \pi^{\frac{n-1}{2}} c^{n+d}}{2d} \frac{\Gamma \left(-\frac{1}{2}-\frac{n}{2d}\right)\Gamma \left(\frac{n}{2d}\right)}{\Gamma \left(\frac{n}{2}\right)} =\mathcal{O}(1) \qquad \textnormal{as } s\to 0_+ \,. \label{eq: d_even_nkd}
\end{align}

If  $ n>d$ then the pole is also simple and one obtains
\begin{align}
G\sim-\frac{\Gamma \left (\frac{1}{2} \right)\Gamma \left (  \frac{n-d}{2} \right)}{\Gamma \left ( \frac{d}{2} \right)}\left( \frac{cs}{2}\right)^{d} \qquad \textnormal{as } s\to 0_+  
\end{align}
and
\begin{align}
\hat\phi (s)\sim 2^{n-1-d}\pi^{\frac{n-1}{2}}c^{2d} \frac{\Gamma \left (\frac{1}{2} \right)\Gamma \left (  \frac{n-d}{2} \right)}{\Gamma \left ( \frac{d}{2} \right)}s^{-n+d} =\mathcal{O}(s^{-n+d}) \qquad \textnormal{as } s\to 0_+ \,.\label{eq: d_even_ngd}
\end{align}

If  $n=d$ then the pole is double and the residue is given by
\begin{align}
G\sim -  \lim_{t\to\frac{1}{2}}    \frac{\partial}{\partial t}\left(\left( t-\frac{1}{2}\right )^2 \frac{    \Gamma \left(t\right)\Gamma \left ( - \frac{1}{2}-t \right) \Gamma \left ( \frac{n}{2}-dt \right) }{\Gamma\left(dt \right)}\left(\frac{cs}{2}\right)^{2dt} \right) \qquad \textnormal{as } s\to 0_+ .
\end{align}
Straightforward calculations lead to
\begin{align}
G\sim \frac{ 2^{-n} \sqrt{\pi}(c s)^{n}  }{d\Gamma\left(\frac{n}{2}\right)}      \Biggl( d\gamma-1 +2d\log (cs)-d\Psi^0 \left(\frac{n}{2}\right)+-2\log\left(2\right) \Biggr) \qquad \textnormal{as } s\to 0_+ 
\end{align}
and
\begin{align}
\hat\phi (s)&\sim  - \frac{\pi^{\frac{n}{2}} c^{n+d}  }{ 2d \Gamma\left(\frac{n}{2}\right)}   \Biggl( d\gamma-1 +2d\log (cs)-d\Psi^0 \left(\frac{n}{2}\right)+-2\log\left(2\right) \Biggr) \qquad \textnormal{as } s\to 0_+ \\ 
&=\mathcal{O}\left(-\log(s)\right) \,  \label{eq: d_even_ngleichd},
\end{align}
where $\Psi^0$ denotes the Digamma function and $\gamma$ is the Euler–Mascheroni constant.

To reproduce polynomials of orders higher than zero, it is sufficient that the RBF's Fourier transform has a singularity at zero of at least order $n+1$. Moreover, we only consider singularities with even parity.  Based on the results from Theorem \ref{Satz: asymptoticbehaviour} this is only possible when $d$ is odd. 
 For greater $d$ we get higher order singularities, so we hope to get also a higher polynom reproduction. 
Now we know that $d$ should be not even, we will only focus the case where $d$ is odd, because we like to have singularities of integer orders. Like in the series expansion of the generalized Fourier transform of the multiquadric there can appear logarithmic terms. So we need to know if and if so at which order does the logarithmic term appears. This will result in a first limitation to the order of polynomial reproduction.

\begin{theorem}
Let $d$ and $n$ be odd. The asymptotic series of $\hat\phi(s)$ for $s\to 0_+$ up to the first logarithmic term is given by
\begin{align}
\hat\phi (s)\sim& \sum_{j=0}^{\left\lceil \frac{n}{2 d}-\frac{1}{2}\right\rceil}C_j(c,n,d)s^{-n+(2j-1)d}\\
&+\sum_{j=0}^{d \left\lceil \frac{n-d}{2 d}\right\rceil +\frac{d-n}{2}} \tilde C_j(c,n,d) s^{2j}\\
&+s^{-n+2 d \left(\left\lceil \frac{n-d}{2 d}\right\rceil +\frac{1}{2}\right)} \log (cs)\qquad \textnormal{as } s\to 0_+\, , \label{eq: asywewrew}
\end{align}
where $C_j$ and $\tilde C_j$ are some constants depending on $c$, $n$, $d$ and $j$.
For even dimensions there appears no logarithmic term and one can write
\begin{align}
\hat\phi (s)= \sum_{j=0}^{\infty}C_j(c,n,d)s^{-n+(2j-1)d}+\sum_{j=0}^{\infty} \tilde C_j(c,n,d) s^{2dj} \,.
\end{align}
\end{theorem}

\subsection*{\textit{Proof}}
Let the dimension $n$ be odd and remember equation~(\ref{eq: G})

\begin{align}
G=\frac{d^{-n/2}}{2\pi \mathrm{i}}\int_L\frac{\Gamma \left ( - \frac{1}{2}-t \right)\Gamma \left (t \right)\Gamma \left (  \frac{n}{2}-dt \right)}{\Gamma \left ( dt \right)}\left( \frac{cs}{2}\right)^{2dt} \mathrm{d}t\, .
\end{align}
 Then the first double pole is at $t=\frac{m}{2}$, where $m$ is a specific odd integer to be determined. The negative argument of the third Gamma function is
\begin{align}
 dt-\frac{n}{2} = \frac{dm-n}{2} \, ,	\label{eq: dfjsf}
\end{align}
where $m, n$, and $d$ are odd integers. The multiplication of two odd numbers is odd, while the subtraction of two odd numbers is even. Dividing by two we obtain an integer. The first higher order pole is at $t_0=\left\lceil \frac{n}{2d}-\frac{1}{2}\right\rceil +\frac{1}{2}$. Summing up the simple poles less than $t_0$ leads to the two sums in equation~(\ref{eq: asywewrew}), where the first sum belongs to all half integers in $ \left[ -\frac{1}{2},t_0\right)$ and the second sum evaluates the poles at $\{\frac{n}{2d},\frac{n+2}{2d},\dots, t_0-\frac{1}{d}\}$. The logarithmic term comes from the double pole at $t_0$.

If the dimension $n$ is even, then there is no double pole because equation~(\ref{eq: dfjsf}) will not evaluate to an integer. Summing up all residues of the poles get the desired result.

\section{Polynom reproduction}

\begin{theorem}
Let $n$ and $d$ be positive odd integers. Then there exist coefficients $\mu_k$  such that the quasi-interpolant $$Q_hf(x)=\sum_{j\in\mathbb{Z}^n}f(jh)\Psi(x/h-j)$$ using the quasi-Lagrange function $$\Psi(x)=\sum_{k\in\mathbb{Z}^n}\mu_k \sqrt{c^{2d}+||x-k||^{2d}}$$ reproduces polynomials of degree $2d-1$. 
\end{theorem}

\subsection*{\textit{Proof}}
Following the theory described in \cite{qi} one uses $$\Psi(x)=\sum_{k\in\mathbb{Z}^n}\mu_k \phi(||x-k||)\,,$$  where the sum can be chosen finite. With the requirement that we get at least a partition of unity this is possible for odd dimensions $n$.
The Fourier transform is given by
\begin{align}
\hat\Psi(y)=\sum_{k\in\mathbb{Z}^n} \mu_k \textnormal{e}^{-\mathrm{i}k\cdot y}\hat\phi(y)\,.
\end{align}
The first limitation of the degree of polynomial reproduction is given by $2 d \left(\left\lceil \frac{n-d}{2 d}\right\rceil+1\right)-1$ which is the difference of exponent of the leading order and the logarithmic prefactor minus one in equation~(\ref{eq: asywewrew}).
The coefficients in the  trigonometric polynomial $\sum_{k\in\mathbb{Z}^n} \mu_k \textnormal{e}^{-ik\cdot y}$ have to be chosen such that the first nonvanishing term in Taylor expansion is of order $n+d$. This condition is called moment condition, and it holds by definition
 if and only if $\sum_{k\in\mathbb{Z}^n} \mu_k k^\alpha=0,  \forall \alpha\in\mathbb{Z}^n , |\alpha|<n+d$.  Then Condition~2 of the Strang and Fix theorem is true. Furthermore, the polynomial is $2\pi $-periodic and the Fourier transform $\hat\phi (x)$ evaluated away from zero is non singular, so Condition~3 is also true. Condition 4 can be fulfilled for the first $2 d \left(\left\lceil \frac{n-d}{2 d}\right\rceil+1\right)-1$ terms, if the coefficients are chosen correctly. The technically challenging part is to show that Condition~1 holds. To show the asymptotic behaviour of $|\Psi(x)|$ we need two identities:
\begin{enumerate}
\item The generalized binomial theorem $(x+y)^r =\sum_{k=0}^\infty  \left({r\atop k}\right) x^{r-k}y^k$, where $|x|>|y|$, and $x, y \in \mathbb{R}.$
\item For $x,a \in \mathbb{R}$ and $x>a$ we have $\frac{1}{x+a}=\frac{1}{x}\sum_{k=0}^\infty \left(-\frac{a}{x}\right)^k$.
\end{enumerate}

We rewrite or quasi-Lagrange function as
\begin{align}
\Psi(x)&=\sum_{k\in\mathbb{Z}^n}\mu_k \sqrt{c^{2d}+||x-k||^{2d}}\\
&=\sum_{k\in\mathbb{Z}^n}\mu_k  ||x-k||^d\sqrt{\left(\frac{c}{||x-k||}\right)^{2d}+1}.
\end{align}
Using the first identity and changing the order of summation which we may do because we assume the coefficient $\mu_k$ to be compactly supported,  we get

\begin{align}
\Psi(x)&=\sum_{j=0}^\infty \sum_{k\in\mathbb{Z}^n}\mu_k ||x-k||^d  \left({1/2\atop j}\right ) \left(\frac{c}{||x-k||}\right)^{2dj}\\
&=\underbrace{ \sum_{k\in\mathbb{Z}^n}\mu_k ||x-k||^d}_{=(A)}+\underbrace{ \sum_{j=1}^\infty \sum_{k\in\mathbb{Z}^n}\mu_k  \left({1/2\atop j}\right ) \frac{c^{2dj}}{||x-k||^{(2j-1)d}}}_{=(B)}\label{eq: dfsagsgh}.
\end{align}
The first part $(A)$ is the well known polyharmonic spline. The decay rate can be arbitrarily fast, depending on the choice of the trigonometric polynomial \cite{BUHMANN2015156}. More precisely, the decay rate is limited by the first nonvanishing term withorder greater than the given singularity of the RBF at zero. In our case equation~(\ref{eq: asywewrew}) shows that this order should be $\min \{n+3d , 2n+2d\}$. Hence there exist coefficients $\mu_k$ such that the polyharmonic spline decays like $\mathcal{O}\left(||x||^{\min \{n+3d-1 , 2n+2d-1\}}\right)$.

To analyse the decay rate of $(B)$  from equation~(\ref{eq: dfsagsgh}) let $x_i$ be the largest component of $x$. 
If $||x||$ goes to infinity, $|x_i|$  will go to infinity as well. Writing this as an asymptotic expansion means that there exists a constant $C_x \in [1,\sqrt{n}]$ depending only on the direction how $x$ tends to infinity, such that $||x|| \sim C_x |x_i|$. Therefore we can write

\begin{align}
(B)&\sim C_x \sum_{j=1}^\infty c^{2dj} \left({1/2\atop j}\right ) \sum_{k\in\mathbb{Z}^n}\mu_k \displaystyle \operatorname {sgn}(x_i-k_i)  \left(\frac{1}{x_i-k_i}\right)^{(2j-1)d} \qquad \textnormal{as } ||x|| \to \infty\\
&=C_x \sum_{j=1}^\infty c^{2dj} \left({1/2\atop j}\right ) \sum_{k\in\mathbb{Z}^n}\mu_k \displaystyle \operatorname {sgn}(x_i-k_i)   \left( \frac{1}{x_i} \sum_{\ell=0}^\infty \left(\frac{k_i}{x_i}\right)^\ell   \right)^{(2j-1)d}\\
&= C_x\sum_{j=1}^\infty c^{2dj} \left({1/2\atop j}\right ) \sum_{k\in\mathbb{Z}^n}\mu_k\displaystyle \operatorname {sgn}(x_i-k_i) \frac{1}{x_i^{(2j-1)d}}   \left(  \sum_{\ell=0}^\infty \left(\frac{k_i}{x_i}\right)^\ell   \right)^{(2j-1)d}\\
&= C_x\sum_{j=1}^\infty c^{2dj} \left({1/2\atop j}\right ) \sum_{k\in\mathbb{Z}^n}\mu_k\displaystyle \operatorname {sgn}(x_i-k_i)  \frac{1}{x_i^{(2j-1)d}}  \sum_{\ell=0}^\infty \sum_{{0\leq s_1,\dots, s_{(2j-1)d}\leq \ell\atop s_1+\dots+s_{(2j-1)d}=\ell}}\left(\frac{k_i}{x_i}\right)^\ell\\
&=C_x \sum_{j=1}^\infty c^{2dj} \left({1/2\atop j}\right )  \sum_{{0\leq s_1,\dots, s_{(2j-1)d}\leq \ell\atop s_1+\dots+s_{(2j-1)d}=\ell}} \frac{1}{x_i^{(2j-1)d}} \sum_{\ell=0}^\infty  \frac{1}{x_i^\ell} \sum_{k\in\mathbb{Z}^n} \displaystyle \operatorname {sgn}(x_i-k_i) \mu_k k_i^\ell\\
&=\mathcal{O}(x_i^{-2d-n})=\mathcal{O}(||x||^{-2d-n}).
\end{align}
Here, the upper bound on $C_x \leq\sqrt{n}$ is essential.
Since the $\mu_k$ are compactly supported with respect to $k$, $x_i$ will become larger than every $k_i$, so that we can use the  moment conditions $ \sum_{k\in\mathbb{Z}^n} \mu_k k^\ell=0$ for all $\ell<n+d$. The lowest nonvanishing term is given for $j=1$ and $\ell=n+d$. 
The overall decay rate is then limited by $(B)$ and is given by $\Psi(x)=\mathcal{O}(||x||^{-2d-n})$ as $||x|| \to \infty$. The second limitation of the degree of polynomial reproduction is given by $2d-1$ which ensures that the sum of the quasi-interpolant with a polynomial of degree $2d-1$ converges. 
Summing up, the polynomial reproduction is limited by the decay rate of $\Psi(x)$ and is given by $2d-1$. Note that the decay rate can be further improved by using linear combinations of our quasi-Lagrange function \cite{qi} and so even higher polynomial reproduction is conceivable.

\section{Error estimates}
Using the Strang and Fix conditions in  Theorem \ref{SF}, we showed that the error estimate is given by 
\begin{align}
||Q_hf-f||_\infty=\mathcal{O}(h^{2d}\log(1/h)) \qquad \textnormal{as } h \to 0.
\end{align}
In contrast the error estimate of the classical multiquadric is given by 
\begin{align}
||Q_hf-f||_\infty=\mathcal{O}(h^{n+1}\log(1/h)) \qquad \textnormal{as } h \to 0.
\end{align}
Comparing these two error estimates reveals improvement when $d>\frac{n+1}{2}$. Furthermore the results for $d=1$ are not equal because our error estimate does not depend on the dimension $n$. These findings indicate that although we have made an improvement, there is still room for further improvement.

\section{Summary}
Summarizing the results, we have discovered an integral representation for the Fourier transform of our new generalized multiquadrics function $\phi(x)=\sqrt{c^{2d}+||x||^{2d}}$.  Furthermore we presented a Meijer G-function representation of the Fourier transform for integer $d$ as a special but useful case. Although not required for quasi-interpolation, we showed that the Fourier transform is not positive definite.
For odd dimensions $n$ and odd $d$ we showed
\begin{align}
\hat\phi (x)\sim -2^{n+d} \pi ^{\frac{n}{2}} \frac{\Gamma \left(\frac{d+n}{2}\right)}{\Gamma \left(-\frac{d}{2}\right)} ||x||^{-n-d} \qquad \textnormal{as } x\to 0 
\end{align}
Furthermore the asymptotic behaviour of the quasi-Lagrange function $\Psi(x)=\sum_{k\in\mathbb{Z}^n}\mu_k\phi (||x-k||)$ for large argument was determined by $$|\Psi(x)|=\mathcal{O}(||x||^{-2d-n}) \qquad \textnormal{as } ||x|| \to \infty.$$
 With these results, we applied the Strang and Fix conditions and found that using this new generalization of multiquadrics as RBF we can construct a quasi-Lagrange function that can reproduce all polynomials of degree $2d-1$. Our quasi-interpolant satisfies the error estimate
\begin{align}
||Q_hf-f||_\infty=\mathcal{O}(h^{2d}\log(1/h)) \qquad \textnormal{as } h \to 0.
\end{align}

\section{Conclusion and Conjecture}
We showed that  in quasi-interpolation the generalized multiquadric can be used, like the classical multiquadric, in odd dimensions. We introduced an arbitrary parameter $d\in \mathbb{N}$, such that the error estimate is improved when $d>\frac{n+1}{2}$. 
We conjecture that this result can be further improved to an error estimate of $\mathcal{O}(h^{n-1+2d}\log(1/h)) $.

The presented theory is only usable in odd dimensions, therefore is would be of great interest to have a second variant which can be used in even dimensions. One candidate would be a generalized version of the shifted thin plate spline $\phi(x)=(c^{2d}+||x||^{2d})\log(c^{2d}+||x||^{2d})$ and will be part of further research.


%
%
%
%
%
%
%
%

\bibliographystyle{abbrv}
\bibliography{Quellen}

\end{document}